\documentclass[12pt]{amsart}

\newtheorem{thm}{Theorem}
\newtheorem{defin}{Definition}

\newtheorem{corol}{Corollary}
\newtheorem{rem}{Remark}
\newtheorem{exa}{Example}

\newcommand{\N}{{\mathbb N}}

\newcommand{\R}{{\mathbb R}}

\newcommand{\car}{\text{char}} 

\newcommand{\iso}{\text{Iso}}   

\newcommand{\spec}{\text{spec}} 

\newcommand{\tr}{\text{tr}\,}

\newcommand{\tw}{\tilde{w}}

\newcommand{\Ga}{\Gamma}

\newcommand{\la}{\lambda}

\begin{document}

\bibliographystyle{plain}

\title[Matrices and Determinants]
{A family of determinants associated with a square matrix}

\author{Eugene Gutkin}

\address{UMK and IMPAN \\
Chopina 12/18\\
87 -- 100 Torun\\
Poland}

\email{\ gutkin@mat.uni.torun.pl}


\date{\today}

\begin{abstract}
We associate with a matrix over an arbitrary field an infinite
family of matrices whose sizes vary from one to infinity; their
entries are traces of powers of the original matrix. We explicitly
evaluate the determinants of matrices in our family. The work is
motivated by applications to graph spectra.
\end{abstract}

\maketitle

\section{A two-parameter family of matrices}       \label{intro}
Let $K$ be a field, and let $G$ be a square matrix of arbitrary
finite size, with entries in $K$. We will denote by $\tr A$ the
trace of a square matrix. Let $t\ge 1$ and $l\ge 0$ be any
integers. Set
\begin{equation}   \label{matrices_eq}
M_{t,l}=M_{t,l}(G)=\left(
\begin{matrix}
\tr G^l & \tr G^{l+1} & \cdots & \tr G^{l+t-1} \\
\tr G^{l+1} & \tr G^{l+2} & \cdots & \tr G^{l+t} \\
\cdot & \cdot &  \cdots & \cdot \\
\tr G^{l+t-1} & \tr G^{l+t} & \cdots & \tr G^{l+2t-2}
\end{matrix}
\right).
\end{equation}

If $A$ is a matrix, we will denote by $[A]_i^j$ the matrix entry
in the $i$th row and $j$th column. In this notation we have
\begin{equation}   \label{entry_eq}
[M_{t,l}(G)]_i^j=\tr G^{i+j+l-2}:\,1 \le i,j \le t.
\end{equation}

\medskip

Matrices similar to these arise in the study of graph spectra
\cite{Ch97,Gut08}. The purpose of this note is to evaluate the
determinants $\det M_{t,l}$.


\vspace{4mm}

\medskip

\section{Evaluation of determinants}   \label{evalu}
We will denote by $\bar{K}$ the algebraic closure of $K$. Let
$t\ge 1$, and let $z_1,\dots,z_t\in \bar{K}$ be arbitrary numbers.
We denote by $V(z_1,\dots,z_t)$ the {\em Vandermonde matrix}:
$$
V(z_1,\dots,z_t)=\left(
\begin{matrix}
1 & 1 & \cdots & 1 \\
z_1 & z_2 & \cdots & z_t \\
\cdot & \cdot &  \cdots & \cdot \\
z_1^{t-2} & z_2^{t-2} & \cdots & z_t^{t-2} \\
z_1^{t-1} & z_2^{t-1}  & \cdots & z_t^{t-1}
\end{matrix}
\right).
$$
%

\medskip

Let $A:K^n\to K^n$ be a linear map. By the {\em multiplicity} of
its eigenvalue $\la\in\bar{K}$ we will always mean the {\em
algebraic multiplicity}.

\medskip
\begin{thm}    \label{det_thm}
Let $G:K^n\to K^n$ be a square matrix of any size. Let
$\la_1,\dots,\la_m\in\bar{K}$ be the distinct eigenvalues of $G$.
Let $p_1,\dots,p_m\in\N$ be their respective multiplicities.

\noindent 1. If $t>m$ then $\det M_{t,l}=0$.

\noindent 2. Let $T\subset\{1,\dots,m\}$ be a nonempty subset and
let $0<t\le m$ be its cardinality. If $T=\{i_1<\cdots<i_t\}$, we
set
$$
p(T)=p_{i_1}\cdots p_{i_t},\ \la(T)=\la_{i_1}\cdots\la_{i_t}.
$$

Then for $t\le m$ we have
\begin{equation}    \label{det_eq}
\det
M_{t,l}=\sum_{T\subset\{1,\dots,m\},|T|=t}p(T)(\la(T))^l{\det}^2
V(\la_{i_1},\dots,\la_{i_t}).
\end{equation}
\begin{proof}
Let $\la$ be an arbitrary number, and let $t>0$. Set
\begin{equation}    \label{lambda_eq}
\vec{\la}=\left[
\begin{matrix}
1\\
\la\\
\la^2\\
\cdot\\
\cdot\\
\cdot\\
\la^{t-1}
\end{matrix}
\right].
\end{equation}
Thus, if $\la\in\bar{K}$, then $\vec{\la}\in\bar{K}^t$. For any
scalar $c$ the notation $c\vec{\la}$ has the usual meaning. We
will in particular use this notation when $c=\la^s,\,0\le s,$
yielding vectors $\la^s\,\vec{\la}$.

We will often view a $t\times t$ matrix as a collection of $t$
column vectors, and use the corresponding notation. For instance,
$$
V(\la_1,\dots,\la_t)=\left(\vec{\la}_1,\dots,\vec{\la}_t\right),\
\det(V(\la_1,\dots,\la_t))=\det\left(\vec{\la}_1,\dots,\vec{\la}_t\right).
$$

For any $k\ge 0$, we have
$$
\tr G^k = \sum_{i=1}^m p_i\la_i^k.
$$
Hence, in our notation,
\begin{equation}    \label{matrix_eq}
M_{t,l}=\left(\sum_{i=1}^mp_i\la_i^{l}\vec{\la}_i,\sum_{i=1}^mp_i\la_i^{l+1}\vec{\la}_i,\dots,\sum_{i=1}^mp_i\la_i^{t+l-1}\vec{\la}_i\right).
\end{equation}

We expand $\det M_{t,l}$ over the columns. Let $i_1,\dots,i_t$ be
$t$ arbitrary indices between $1$ and $m$. We denote by
$\vec{i}=(i_1,\dots,i_t)$ the corresponding multi-index. Then, by
equation~\eqref{matrix_eq}
$$
\det M_{t,l} = \sum_{(i_1,\dots,i_t)}\det\left(
p_{i_1}\la_{i_1}^{l}\vec{\la}_{i_1},p_{i_2}\la_{i_2}^{l+1}\vec{\la}_{i_2},\dots,p_{i_t}\la_{i_t}^{l+t-1}\vec{\la}_{i_t}\right).
$$
If not all indices $i_1,\dots,i_t$ are distinct, the matrix
$\left(\vec{\la}_{i_1},\vec{\la}_{i_2},\dots,\vec{\la}_{i_t}\right)$
is degenerate; thus, the corresponding contribution to $\det
M_{t,l}$ is zero. When $m<t$, then the indices $i_1,\dots,i_t$
cannot be all distinct. This proves our first claim.

\medskip

From now on $1\le t \le m$. Let $I=\{(i_1,\dots,i_t)\}$ be the set
of distinct multi-indices. The preceding observation yields
\begin{equation}   \label{expansion_eq}
\det M_{t,l} = \sum_{(i_1,\dots,i_t)\in I}p_{i_1}\cdots
p_{i_t}(\la_{i_1}^l\la_{i_2}^{l+1}\cdots\la_{i_t}^{l+t-1})\det\left(\vec{\la}_{i_1},\vec{\la}_{i_2},\dots,\vec{\la}_{i_t}\right).
\end{equation}
Every subset $T\subset\{1,\dots,m\}$ of cardinality $t$ determines
$t!$ multi-indices in $I$. Precisely one of them goes in the
increasing order: $i_1(T)<\cdots<i_t(T)$; we will denote this
multi-index by $\vec{i}(T)$. On the other hand, any multi-index
$\vec{i}=(i_1,\dots,i_t)\in I$ determines a subset
$T=T(\vec{i})=\{i_1,\dots,i_t\}$ of cardinality $t$ and a
one-to-one mapping $w=w(\vec{i}):T\to T$ defined by
\begin{equation}   \label{permut_eq}
w:i_1(T)\mapsto i_1,\dots,i_t(T)\mapsto i_t.
\end{equation}

\medskip

We will denote by $(-1)^w$ the usual sign function on
permutations. Thus, $(-1)^w=1$ (resp. $(-1)^w=-1$) if the
permutation $w$ is even (resp. odd). Let
$\vec{i}=(i_1,\dots,i_t)\in I$ and let
$T=T(\vec{i}),\,w=w(\vec{i})$. Then $p_{i_1}\cdots p_{i_t}$ and
$\la_{i_1}^l\la_{i_2}^l\cdots\la_{i_t}^l$ depend only on $T$. We
have $p_{i_1}\cdots
p_{i_t}=p(T),\,\la_{i_1}^l\la_{i_2}^l\cdots\la_{i_t}^l=(\la(T))^l$.
Set
$v(T)=\det\left(\vec{\la}_{i_1(T)},\vec{\la}_{i_2(T)},\dots,\vec{\la}_{i_t(T)}\right)$.
Then
$$
\det\left(\vec{\la}_{i_1},\vec{\la}_{i_2},\dots,\vec{\la}_{i_t}\right)=(-1)^{w(\vec{i})}v(T).
$$

Using these relationships, we rewrite the above expansion for
$\det M_{t,l}$ as a repeated sum: the former summation is over the
sets $T\subset\{1,\dots,m\}$ of cardinality $t$, while the latter
is over the $t!$ permutations of $T$. Moreover, we identify the
set $\iso(T)$ of permutations of elements in $T$ with the
symmetric group $S_t$ as follows. Let $i_1(T)<\cdots<i_t(T)$ be
the elements of $T$ in the increasing order, and let $i\in T$ be
arbitrary. Then there is a unique $1\le s=o(i)\le t$ such that
$i=i_s(T)$. We say that $s$ is the {\em order of $i$ in $T$}. The
order function $i\mapsto o(i)$ identifies $T$ with
$\{1,\dots,t\}$, and hence $\iso(T)$ with $S_t$. Let $w\in\iso(T)$
be given by equation~\eqref{permut_eq} and let $\tw\in S_t$ be the
corresponding permutation. Then
$$
\tw:1\mapsto o(i_1),\dots,\tw:t\mapsto o(i_t).
$$

Equation~\eqref{expansion_eq} yields
\begin{equation}   \label{doubl_sum_eq}
\det M_{t,l} =
\sum_{T\subset\{1,\dots,m\},|T|=t}p(T)(\la(T))^lv(T)\sum_{w\in
\iso(T)}(-1)^w\la_{i_1}^0\la_{i_2}^1\cdots\la_{i_t}^{t-1}.
\end{equation}
We recall the well known formula
$$
\det\left(\vec{\la}_1,\dots,\vec{\la}_t\right)=\sum_{w\in S_t
}(-1)^w\la_{w(1)}^{0}\la_{w(2)}^{1}\cdots\la_{w(t)}^{t-1}.
$$
Using this identity and the preceding identification of $S_t$ and
$\iso(T)$, we obtain
$$
\sum_{w\in\iso(T)}(-1)^w\la_{i_1}^0\la_{i_2}^1\cdots\la_{i_t}^{t-1}=\det\left(\vec{\la}_{i_1(T)},\vec{\la}_{i_2(T)},\dots,\vec{\la}_{i_t(T)}\right).
$$

\medskip

Thus, we have evaluated the second sum in
equation~\eqref{doubl_sum_eq}. Now we rewrite the expansion
equation~\eqref{doubl_sum_eq} as
\begin{equation}   \label{answer_eq}
\det M_{t,l} =
\sum_{T\subset\{1,\dots,m\},|T|=t}p(T)(\la(T))^lv^2(T),
\end{equation}
which is a short form for equation~\eqref{det_eq}.
\end{proof}
\end{thm}

\medskip

\medskip

\medskip

\vspace{3mm}

\section{Examples and corollaries}   \label{examp}
In order to illustrate Theorem~\ref{det_thm}, we will now consider
matrices with very small numbers of eigenvalues.
\begin{exa}     \label{small_exa}
{\em

\noindent i) Let $G$ be a matrix with one eigenvalue, say $\la$.
This is the special case $m=1$ in equation~\eqref{det_eq}. Then
the multiplicity, say $p$, of $\la$ coincides with the size of
$G$. We have
\begin{equation}   \label{one_eigen_eq}
M_{t,l}(G)=\left(
\begin{matrix}
p\la^l & p\la^{l+1} & \cdots & p\la^{l+t-1} \\
p\la^{l+1} & p\la^{l+2} & \cdots & p\la^{l+t} \\
\cdot & \cdot &  \cdots & \cdot \\
p\la^{l+t-1} & p\la^{l+t} & \cdots & p\la^{l+2t-2}
\end{matrix}
\right).
\end{equation}
This matrix is nondegenerate only if $t=1$, yielding $\det
M_{1,l}=p\la^l$ and $\det M_{t,l}=0$ for $t>1$.

\medskip

\noindent ii) Let $G$ be a matrix with two eigenvalues, say $\la$
and  $\mu$. We denote by $p$ and $q$ their respective
multiplicities. We have
$$
M_{t,l}(G)=\left(
\begin{matrix}
p\la^l+q\mu^l & p\la^{l+1}+q\mu^{l+1} & \cdots & p\la^{l+t-1}+q\mu^{l+t-1} \\
p\la^{l+1}+q\mu^{l+1} & p\la^{l+2}+q\mu^{l+2} & \cdots & p\la^{l+t}+q\mu^{l+t} \\
\cdot & \cdot &  \cdots & \cdot \\
p\la^{l+t-1}+q\mu^{l+t-1} & p\la^{l+t}+q\mu^{l+t} & \cdots &
p\la^{l+2t-2}+q\mu^{l+2t-2}
\end{matrix}
\right).
$$
By straightforward calculations, $\det M_{1,l}=p\la^l+q\mu^l$ and
$\det M_{2,l}=pq(\la\mu)^l(\la-\mu)^2$. This corresponds to
$t=1,2$  in equation~\eqref{det_eq}. For $t>2$ the above matrix is
degenerate.

}
\end{exa}

\medskip

For applications of Theorem~\ref{det_thm}, we are especially
interested in the case when $G$ is a real valued symmetric matrix.
Then the eigenvalues of $G$ are real. Moreover, the algebraic
multiplicities of eigenvalues of $G$ coincide with their geometric
multiplicities \cite{Ga,HJ}. By $M_t=M_t(G)$ we will mean the
matrix $M_{t,0}(G)$.

\medskip
\begin{corol}   \label{real_cor1}
Let $G$ be a real, symmetric square matrix of an arbitrary size.
Suppose that $G$ has $m$ distinct eigenvalues. Then $\det
M_t(G)>0$ for $t\le m$ and $\det M_t(G)=0$ for $t>m$.
\begin{proof}
Let $\la_1,\dots,\la_m$ be the eigenvalues of $G$. By
Theorem~\ref{det_thm}, for $t\le m$ we have
\begin{equation}   \label{det_eq1}
\det M_t=\sum_{T\subset\{1,\dots,m\},|T|=t}p(T){\det}^2
\left(\vec{\la}_{i_1},\dots,\vec{\la}_{i_t}\right).
\end{equation}
Since $\la_1,\dots,\la_m\in\R$, for any set $T$ in
equation~\eqref{det_eq1} the number $v(T)^2={\det}^2
\left(\vec{\la}_{i_1},\dots,\vec{\la}_{i_t}\right)$ is positive.
Our first claim follows. The other claim is contained in
Theorem~\ref{det_thm}.
\end{proof}
\end{corol}

\medskip

\begin{defin}    \label{spec_size_def}
Let $G$ be a square matrix of any size with entries in an
arbitrary field.  The {\em spectral size} of $G$ is the number of
its distinct eigenvalues.
\end{defin}
\medskip
We will use Theorem~\ref{det_thm} to characterize matrices with a
particular  spectral size.
\begin{corol}   \label{spec_size_cor}
Let $K$ be a field and let $G$ be a square matrix with entries in
$K$. Then the following holds.

\noindent 1. The spectral size of $G$ is equal to $m$ iff $\det
M_m(G)\ne 0$ and $\det M_t(G)=0$ for $t>m$.

\noindent 2. Let the spectral size of $G$ be equal to $m$, and let
$l\ge 1$. Then
$$
\det M_{m,l}(G)=0
$$
iff $G$ is a degenerate matrix.
\begin{proof}
Let $m$ be the spectral size of $G$, and let $s\ge 0$ be any
integer. By Theorem~\ref{det_thm}, $\det M_t(G)=0$ for all $t\ge
s$ iff $s>m$. Also, by Theorem~\ref{det_thm}, $\det M_m(G)\ne 0$.
This proves claim one.

Let $\la_1,\dots,\la_m$ be the eigenvalues of $G$. By
equation~\eqref{det_eq}
\begin{equation}   \label{prod_eigen_eq}
\det M_{m,l}=\left(\la_1\cdots\la_m\right)^l\det M_m.
\end{equation}
We have shown already that $\det M_m\ne 0$. Hence $\det M_{m,l}=0$
iff zero is an eigenvalue of $G$. Claim two follows.
\end{proof}
\end{corol}

\medskip

Let $G$ be a square matrix over a field $K$. Let
$\la_1,\dots,\la_m\in\bar{K}$ be its distinct eigenvalues. Let
$p_1,\dots,p_m$ be their respective multiplicities. The monic
polynomial $P_{\car}(\la)=\prod_{i=1}^m(\la-\la_i)^{p_i}$ over $K$
is the {\em characteristic polynomial} of $G$ \cite{HJ}. We have
$P_{\car}(G)=0.$ The monic polynomial $P_{\min}$ of minimal degree
satisfying $P_{\min}(G)=0$ is the {\em minimal polynomial} of $G$
\cite{Ga}. We set $P_{\spec}(\la)=\prod_{i=1}^m(\la-\la_i)$; we
call $P_{\spec}(\la)$ the {\em spectral polynomial} of
$G$.\footnote{If a danger of confusion arises, we will indicate
the dependence on $G$ by superscripts, e.g.,
$P_{\spec}^{(G)}(\la)$ for the spectral polynomial.}

\medskip

We will use the notation $I$ for the identity matrix of any size.
The size of a particular identity matrix should be clear from the
context.

\begin{corol}   \label{character_cor}
Let $G$ be an arbitrary square matrix over a field $K$. Then
\begin{equation}   \label{mini_eq}
P_{\spec}(\la)=\frac{\det M_{m,1}(\la I-G)}{\det M_m(G)}.
\end{equation}
If $G$ is a symmetric square matrix over $\R$, then the polynomial
in the right hand side of equation~\eqref{mini_eq} is the minimal
polynomial of $G$.
\begin{proof}
Consider the matrices $\la I-G,\,\la\in K$. The eigenvalues of
$\la I-G$ are $\la-\la_1,\dots,\la-\la_m$; their multiplicities do
not depend on $\la$. Thus, by Theorem~\ref{det_thm}, for all $t$
we have the identities
$$
\det M_t(\la I-G) = \det M_t(G).
$$
Our first claim now follows from equation~\eqref{det_eq}.

By definition, the minimal polynomial of $G$ divides the
characteristic polynomial, and the spectral polynomial divides the
minimal polynomial of $G$. Moreover, the minimal polynomial
coincides with the spectral polynomial iff the matrix $G$ is
semi-simple \cite{Ga,HJ}. Since real symmetric matrices are
semi-simple, our second claim follows from the first.
\end{proof}
\end{corol}

\medskip

\begin{rem}    \label{homothet_rem}
{\em

Let $c\in K$, $c\ne 0$, and set $G_1=cG$. Then the eigenvalues and
their multiplicities satisfy
$\la_i(G_1)=c\la_i(G),p_i(G_1)=p_i(G)\,:1 \le i \le m$. It is then
immediate from equation~\eqref{det_eq} that
$$
\det M_{t,l}(cG) = c^{ml+t(t-1)} \det M_{t,l}(G).
$$
}
\end{rem}
%


\medskip


%

%
\medskip


\section{Concluding remarks and amplifications}
Let $G_1,G_2$ be the adjacency matrices of finite rooted graphs.
Let $\Ga$ be their {\em free product}.\footnote{This material
extends to the free products of any number of graphs.} The Green
function of $\Ga$ can be expressed in terms of the spectra of
$G_1,G_2$ \cite{Gut98}. This expression yields some information
about the spectrum of $\Ga$ \cite{Gut98}. Complete information
about the spectrum of $\Ga$ can be obtained this way if $G_1,G_2$
have sufficiently small spectra. The present work will be used to
analyze free products of graphs with small numbers of eigenvalues
\cite{Gut08}.

\medskip

Besides applications to graph spectra, we feel that the matrices
in equation~\eqref{matrices_eq} are of interest on their own. In
particular, it seems natural to extend the preceding material to
arbitrary reflection groups. See \cite{Bou,Hum} and \cite{Gut73}.
The analysis of matrices in equation~\eqref{matrices_eq} is based,
in a certain sense, on the symmetric group, which is the standard
example of a reflection group. Note that there are other examples
of matrix families naturally associated with reflection groups
\cite{Gut73}.

\vspace{3mm}

\noindent{\bf Acknowledgements}: Some of the present work was
performed in July-September 2008, while the author was visiting
UCLA. It is a pleasure to thank the UCLA mathematical department
for the hospitality.

\medskip




\end{document}